\documentclass{amsart}
\usepackage{}
\usepackage{epsfig}
\usepackage{amsmath}
\usepackage{amsfonts}
\begin{document}

\newtheorem{thm}{Theorem}[section]
\newtheorem{lem}{Lemma}[section]
\newtheorem{cor}{Corollary}[section]
\newtheorem{conj}{Conjecture}
\newtheorem{qn}{Question}
\newtheorem{pro}{Proposition}[section]

 \theoremstyle{definition}
\newtheorem{defn}{Definition}[section]

\theoremstyle{remark}
\newtheorem{rmk}{Remark}[section]

\def\square{\hfill${\vcenter{\vbox{\hrule height.4pt \hbox{\vrule width.4pt
height7pt \kern7pt \vrule width.4pt} \hrule height.4pt}}}$}

\def\Z{\Bbb Z}
\def\R{\Bbb R}
\newenvironment{pf}{{\it Proof:}\quad}{\square \vskip 12pt}

\title[complex hyperbolic manifolds]{Balls in  complex hyperbolic manifolds}
\thanks{*Corresponding author}
\author{ Baohua Xie, Jieyan Wang and Yueping Jiang*}

\address{College of Mathematics and Econometrics \\ Hunan University \\ Changsha, 410082, China}
\email{xiexbh@gmail.com, jywang@hnu.edu.cn, ypjiang731@163.com}
 \maketitle

\begin{abstract}
In this paper we get an explicit lower bound for the radius of a Bergman ball contained in the Dirichlet fundamental polyhedron of
a torsion-free discrete group $G\subset\mathbf{PU}(n,1) $ acting on complex hyperbolic space.
As an application, we also give a lower bound for the volumes of complex hyperbolic $n$-manifolds.

 \end{abstract}

\section{Introduction}
Recall that any (real or complex) hyperbolic $n$-manifold can be identified with $\mathbf{H}^{n}/G$, where $G \subset \rm{Isom}(\mathbf{H}^{n})$ is a torsion-free discrete isometric subgroup acting on real or complex hyperbolic $n$-space $\mathbf{H}^{n}$. Let $\mathcal {F}_{n}$ be the set of all discrete torsion-free groups and $o$ be any point in hyperbolic space. Define
$$r_{n}=\frac{1}{2}\inf_{G_{n}\in\mathcal {F}_{n}}\sup_{o\in \mathbf{H}^{n}}\inf_{f\in G_{n}}d_{\mathbf{H}^{n}}\big(o,f(o)\big)$$ where $d_{\mathbf{H}^{n}}$ is hyperbolic distance in the case of real hyperbolic manifolds and is Bergman distance in the case of complex hyperbolic manifolds. It is known that $r_{n}$ is the largest number such that every real or complex hyperbolic $n$-manifold contains an embedded ball of that radius. A positive lower
bound on $r_{n}$ provides geometric informations about all hyperbolic $n$-manifolds, such as the lower bounds for the volumes of hyperbolic manifolds and the thick and thin decomposition of such
manifolds. The famous work of Kazdan and Margulis \cite{km} implied the existence of a uniform hyperbolic ball of radius $r_{n}$, which is embedded in any hyperbolic $n$-manifold, see also \cite{wa}. However, they did not provide a computable value of $r_{n}$.

 For the real hyperbolic $3$-manifolds, the bound $r_{3}\geq\frac{1}{25}$ was produced by  Waterman \cite{wat}. A better result was included in  \cite{gm}.
 In the case of high dimension, Martin \cite{ma1} found an explicit lower bound for $r_{n}$ depending only on the dimension $n$. He obtained this bound by using the generalized J{\o}rgensen's inequality which he gave in paper \cite{ma2}. Then he found a lower bound on the radius of a hyperbolic ball contained in the Dirichlet fundamental polyhedron of the uniformizing group acting on  hyperbolic space.
 Friedland and Hersonsky \cite{fh} improved the constant in Martin's inequality slightly. They also use their results to estimate the radii of balls in hyperbolic $n$-manifolds.
 If $G \subset \rm{Isom}(\mathbf{H}^{n})$ contains torsion elements, $\mathbf{H}^{n}/G$ is a hyperbolic orbifold. For this case, Adeboye \cite{ad} derived an explicit lower bound for the volume of a complete hyperbolic orbifold, dependent on the dimension and the maximum order of torsion in the orbifold's fundamental group.

For the complex hyperbolic setting, in \cite{hp}, Hersonsky and Paulin proved that the smallest volume of a compact complex hyperbolic $2$-manifold is $8\pi^{2}$. Furthermore, they gave a lower bound for the volume of a cusped complex hyperbolic $n$-manifold. Parker proved that the smallest volume of a cusped complex hyperbolic $2$-manifold is $\frac{8\pi^{2}}{3}$ in \cite{pa2}. For the cusped complex hyperbolic manifolds case, these results were improved by different authors, we refer the readers to \cite{ki} and \cite{hw}, for example. For the complex hyperbolic orbifolds, Parker found a lower bound for the radius of the embedded ball in a cusped complex hyperbolic orbifold in \cite{pa}. Recently, X. Fu et.al\cite{flw} generalized
Adeboye's methods and gave a formula to estimate the volume of complex hyperbolic orbifolds. This formula also depends on the dimension and the maximal order of elliptic elements in the fundamental groups of complex hyperbolic orbifolds.

Motivated by the idea of  Martin\cite{ma1}, we study the embedded balls in  complex hyperbolic $n$-manifolds. We get that every complex $n$-manifold contains a Bergman ball with radius $r_{n}=\frac{0.01}{17^{n-1}}$. As an application, we obtain a lower bound for the volume of a complex hyperbolic $n$-manifold.

\section{Complex hyperbolic space}
First, we recall some background on complex hyperbolic geometry. More details can be found in \cite{go}. In this paper we work in the unit ball model.
\subsection{Ball model}
Let $\mathbf{C}^{n,1}$ be the $n+1$-dimensional complex vector space consisting of $n+1$-tuples
$$\mathbf{z}=(z_{1},z_{2},\ldots,z_{n},z_{1+n})$$ with the Hermitian pairing
$$\langle \mathbf{z},\mathbf{w}\rangle= \mathbf{z}J\mathbf{w}^{*}=z_{1}\overline{w_{1}}+z_{2}\overline{w_{2}}+\ldots+z_{n}\overline{w_{n}}-z_{1+n}\overline{w_{1+n}},$$
where $"*"$ denotes the complex Hermitian transpose, $\mathbf{z},\mathbf{w}$ are the column vectors in $\mathbf{C}^{n+1}$ and
$$J= \left(\begin{array}{cc}I_{n}& 0\\
0& -1\end{array}\right).
$$

We consider the subspaces  of  $\mathbf{C}^{n,1}$
$$V_{-}=\{\mathbf{z}\in \mathbf{C}^{n,1}: \langle \mathbf{z}, \mathbf{z}\rangle<0 \},$$
$$V_{0}=\{\mathbf{z}\in \mathbf{C}^{n,1}: \langle \mathbf{z}, \mathbf{z}\rangle=0 \}.$$

Let $\mathbf{P}: \mathbf{C}^{n+1}-\{0\}\longrightarrow \mathbf{C}\mathbf{P}^{n}$ be the canonical
projection. Then complex hyperbolic space is defined to be $\mathbf{H}_{\mathbf{C}}^{n}=\mathbf{P}(V_{-})$ and $\partial\mathbf{H}_{\mathbf{C}}^{n}=\mathbf{P}(V_{0})$ is its boundary.

We define the ball model of complex hyperbolic space by taking the section defined by $z_{n+1}=1$ for the above Hermitian form. That is, if we take column vectors $\mathbf{z}=(z_{1},z_{2},\ldots,z_{n},1)$ in $\mathbf{C}^{n,1}$, then consider what it means for $\langle \mathbf{z},\mathbf{z}\rangle$ to be negative.

For the above Hermitian form, we obtain $z\in \mathbf{H}_{\mathbf{C}}^{n}$ provided:
$$\langle \mathbf{z},\mathbf{z}\rangle=z_{1}\overline{z_{1}}+z_{2}\overline{z_{2}}+\ldots+z_{n}\overline{z_{n}}-1<0.$$
In other words,
$$|z_{1}|^{2}+|z_{2}|^{2}+\ldots+|z_{n}|^{2}<1.$$
Thus $z=(z_{1},z_{2},\ldots,z_{n})$ is in the unit ball in $\mathbf{C}^{n}$. This forms the unit ball model of complex hyperbolic space. The boundary
of the unit ball model is the sphere $S^{2n-1}$ given by
$$|z_{1}|^{2}+|z_{2}|^{2}+\ldots+|z_{n}|^{2}=1.$$
So we can identify the complex hyperbolic space $\mathbf{H}_{\mathbf{C}}^{n}$ with the ball $\mathbf{B}^{2n}$.

We say that the standard lift of a point  $z=(z_{1},z_{2},\ldots,z_{n})$  in the unit ball to $\mathbf{C}^{n,1}$ is the column vector
$\mathbf{z}=(z_{1},z_{2},\ldots,z_{n},1)$ of $\mathbf{C}^{n,1}$ whose first $n$ coordinates are those of $z$ and whose last coordinate is 1. The
standard lift of 0 is the column vector $e_{n+1}=(0,0,\ldots,1)$.

The distance function in $\mathbf{H}_{\mathbf{C}}^{n}$ has the following useful algebraic description in terms of the Hermitian structure on $\mathbf{C}^{n,1}$.
Let $x,y\in\mathbf{H}_{\mathbf{C}}^{n} $ be points corresponding to vectors $\mathbf{x},\mathbf{y}\in \mathbf{C}^{n,1}$. Then the distance between them is given by
$$\cosh^{2}\big(\frac{\rho(x,y)}{2}\big)=\frac{\langle \mathbf{x},\mathbf{y}\rangle \langle \mathbf{y},\mathbf{x}\rangle}{\langle\mathbf{ x},\mathbf{x}\rangle \langle \mathbf{y},\mathbf{y}\rangle}.$$
The biholomorphic isometric group of $\mathbf{H}_{\mathbf{C}}^{n}$ is $\mathbf{PU}(n,1)=\mathbf{U}(n,1)/\mathbf{U}(1)$.There exist three kinds of
holomorphic isometries of $\mathbf{H}_{\mathbf{C}}^{n}$.

{\rm(i)} Loxodromic isometries, each of which fixes exactly two points of $\partial\mathbf{H}_{\mathbf{C}}^{n}$.
One of these points is attracting and the other repelling.

{\rm(ii)} Parabolic isometries, each of which fixes exactly one point of $\partial\mathbf{H}_{\mathbf{C}}^{n}$.

{\rm(iii)} Elliptic isometries, each of which fixes at least one point of $\mathbf{H}_{\mathbf{C}}^{n}$.

Any matrix in $\mathbf{PU}(n,1)$ fixing the origin $e_{n+1}$ is the projectivisation of a block diagonal matrix $\mathbf{U}(n)\times\mathbf{U}(1)$ in
$\mathbf{U}(n,1)$. It has the form
$$A= \left(\begin{array}{ccc} A_{1}& 0\\
0& e^{i\theta}\end{array}\right)
$$
where $A_{1}\in \mathbf{U}(n)$ and $e^{i\theta}\in\mathbf{U}(1) $. Projectivising we may assume that $e^{i\theta}=1$.

In order to obtain a lower bound for  the volume of a complex hyperbolic $n$-manifold, we need the following lemma.
See, for example, Lemma 6.18 on page 108 and Corollary A.3 on page 254 of \cite{gr}.
\begin{lem}The complex hyperbolic volume of a geodesic ball of radius $r_{0}$ is
$$Vol(B(r_{0}))=\frac{4^{n}\sigma_{2n-1}}{2n}\sinh^{2n}(\frac{r_{0}}{2}),$$ where $\sigma_{2n-1}=\frac{2\pi^{n}}{(n-1)!}$ is the Euclidean volume of the unit sphere
$S^{2n-1}\subset \mathbf{C}^{n}$.
\end{lem}

\section{Preliminaries}

We identify the group of complex hyperbolic isometries with the Lie group $\mathbf{SU}(n,1) $ via the usual topological isomorphism. We shall use letters
such as $f,g,h$ to denote complex hyperbolic isometries and letters such as $A,B,C$ to denote the corresponding matrices in $\mathbf{SU}(n,1)$. We denote
by $\mathbf{U}(n)$ the maximal unitary subgroup of $\mathbf{SU}(n,1) $. We say $G$ is discrete if identity is isolated in $G$. We denote the orbit $\{g(x): g\in G\}$ of a point $x$ under $G$ by $G(x)$. The limit set $L(G)$ of a discrete group $G$ is then the set of accumulation points of the orbit of an arbitrary point
$x\in \mathbf{B}^{2n}$. We say that a discrete group $G$ is elementary if $|L(G)|<2$, non-elementary otherwise.

For any element $A$ of $\mathbf{SU}(n,1) $ we denote its operator norm as
$$\| A\|=\max\{|Av|: v\in \mathbf{C}^{n+1}, \ and \ |v|=1\}.$$

Let $\sigma(A)$ denote the spectra of $A$, that is, the set of all eigenvalues of $A$. The number
$$r_{\sigma}(A)=\max_{\lambda\in \sigma(A)}|\lambda|$$
is defined to be the spectral radius of $A$. We will use an alternative definition for the operator norm in this paper, that is,
$$\parallel A\parallel=\sqrt{r_{\sigma}(A^{*}A)}.$$

The operator norm of $A\in \mathbf{SU}(n,1)$ is determined by its eigenvalues. Thus operator norm is a conjugate invariant.
\begin{pro} Let $A\in \mathbf{SU}(n,1)$ and $B\in \mathbf{U}(n)$. Then $$ \parallel B AB^{-1}\parallel=\parallel A\parallel.$$
\end{pro}

One recall that orientation preserving M\"{o}bius transformations isomorphic to the group $\mathbf{SO}^{+}(1,n)$.
Martin\cite{ma2} obtained the following numeric version of the Zassenhaus-Kazdan-Margulis lemma.
\begin{lem}
  If $A$ and $B$ are elements of $G\subset\mathbf{SO}^{+}(1,n)$, then
  $$\max \{\parallel A-Id\parallel,\parallel B-Id\parallel\}\geq 2-\sqrt{3}$$
  and
  $$\max \{\parallel A-Id\parallel,\parallel [A,B]-Id\parallel\}\geq 2-\sqrt{3}$$
  unless $A$ and $B$ lie in the same elementary subgroup of $G$.
\end{lem} Here $[A,B]=ABA^{-1}B^{-1}$ is the multiplicative commutator and the norm $\parallel \cdot \parallel$ is the spectral norm. This is the torsion free version of Martin's generalization of J{\o}rgensen inequality.
This norm version of J{\o}rgensen inequality was generalized  to complex hyperbolic space by B. Dai et al \cite{dfn}.  and S. Kamiya \cite{ka}.

\begin{rmk} This result was improved in Theorem 2.6 of Friedland and
Hersonsky \cite{fh}, where the constant $2-\sqrt{3}$ was replaced by $\tau$, defined to be the unique positive
solution to $2\tau(\tau+1)^{2}=1$, and is approximately 0.2971.
\end{rmk}

The following lemma is a variation of Lemma 3.1. Martin\cite{ma1} got a lower bound on the radius of a hyperbolic ball in a hyperbolic $n$-manifold by using this lemma.
\begin{lem}
  Let  $G$ be a discrete non-elementary torsion free subgroup of $\mathbf{SO}^{+}(1,n)$. Then
  there is an $\alpha\in \mathbf{SO}^{+}(1,n)$ such that
  $$\parallel A\parallel \parallel A-Id \parallel\geq (1-\frac{\sqrt{3}}{2})^{\frac{1}{2}}> \frac{1}{2\sqrt{2}}$$
  for all $A\in \alpha G \alpha^{-1}$.
\end{lem}
\begin{rmk}Friedland and Hersonsky \cite{fh} improved the constant in the above inequality. That is,
$$\parallel A\parallel \parallel A-Id \parallel\geq \omega$$ where $\omega$ is the unique positive root
of $2\omega(2\omega^{2}+1)=1$, and is approximately 0.3854. They proved this inequality in the general setting
by considering normed algehras with an involution.
\end{rmk}
As a special case of Theorem 2.16 of Friedland and Hersonsky \cite{fh}, we state their result in the
setting of complex hyperbolic case. It play an important role in the proof  of our main theorem.
\begin{lem}
  Let  $G$ be a discrete non-elementary torsion free subgroup of $\mathbf{SU}(n,1)$. Then
  there is an $\alpha\in \mathbf{SU}(n,1)$ such that
  $$\parallel A\parallel \parallel A-Id \parallel\geq \omega$$
  for all $A\in \alpha G \alpha^{-1}$.
\end{lem}

 \section{The main result and its proof }
 Since $G$ is a discrete torsion-free group acting  on the complex hyperbolic space, the Bergman distance $\rho(o,f(o))$ must be larger than a certain number. Our approach follows the main idea of Martin's proof. We begin with some lemmas.

 \begin{lem}(\textbf{Dirichlet's pigeon-hole principle}) Let $\theta_{i}\in [0,1], i=1,2,\ldots,m$. Then for each $Q\geq1$ there is $q\leq Q^{m}$ and $p_{i}, i=1,2,\ldots,m$, such that
 $$|\theta_{i}-\frac{p_{i}}{q}|\leq \frac{1}{qQ}.$$
 \end{lem}

 As application of the conjugacy to diagonal form
for an element of $\mathbf{U}(n)$ together with the above principle and an eigenvalue calculation yields the following. One can compare the following lemma with Lemma 4.1 in \cite{fh} and Corollary 3.3 in \cite{ma1}.

\begin{lem} Let $A\in \mathbf{U}(n)$. Then for each $Q> 1$ there is $B\in \mathbf{U}(n)$ such that for some $1\leq q\leq Q^{n-1}$, $B^{q}=Id$ and
$$\parallel A-B\parallel\leq \frac{2\pi}{qQ}.$$
 \end{lem}
 \begin{pf}
 Without loss of generality, we may assume that there exist $B_{1}\in \mathbf{U}(n)$ such that
 $$A=B_{1}A_{1}B_{1}^{-1}$$
 where
 $$A_{1}= \left(\begin{array}{cccccc} e^{i2\pi \theta_{1}}& & & & & \\
&e^{i2\pi \theta_{2}} & & & &\\ & &. & & & \\& & &. & &\\& & & &e^{i2\pi \theta_{n-1}} &\\& & & & &1\end{array}\right)$$
 $ \theta_{i}\in [0,1]$.

 Then there is  $$B=B_{1}B_{0}B_{1}^{-1}\in \mathbf{U}(n)$$ where
 $$B_{0}= \left(\begin{array}{cccccc}  e^{i2\pi \frac{p_{1}}{q}}& & & & & \\
&e^{i2\pi \frac{p_{2}}{q}} & & & &\\ & &. & & & \\& & &. & &\\& & & &e^{i2\pi \frac{p_{n-1}}{q}} &\\& & & & &1\end{array}\right)$$
 $ \theta_{i}\in [0,1]$.

 Thus we have the following from the definition of operator norm.
 \begin{align*}
\parallel A-B\parallel=\parallel A_{1}-B_{0}\parallel &=\sqrt{r_{\sigma}((A_{1}-B_{0})^{*}(A_{1}-B_{0}))}\\
&=\max \sqrt{|(e^{i2\pi \theta_{i}}-e^{i2\pi \frac{p_{i}}{q}})^{2}|}\\
&=\max\sqrt{|2-2\cos(2\pi \theta_{i}-2\pi \frac{p_{i}}{q})|}\\
&=\max\sqrt{|4\sin^{2}(\pi \theta_{i}-\pi \frac{p_{i}}{q})|}\\
&= 2\max|\sin(\pi \theta_{i}-\pi \frac{p_{i}}{q})|\\
&\leq 2\max|\pi \theta_{i}-\pi \frac{p_{i}}{q}|\leq \frac{2\pi}{qQ}.
\end{align*}

 \end{pf}

\begin{rmk}The similar lemma in real hyperbolic case  was provided in \cite{fh}.

\end{rmk}

 \begin{lem} Let the complex hyperbolic isometry $f$ correspond to the matrix $A\in \mathbf{SU}(n,1) $. Then
  $$\rm{dist}(A,\mathbf{U}(n))=\inf\{\parallel A-O\parallel, O\in \mathbf{U}(n)\}\leq r(r-1),$$
  where $r=\exp(\rho(0,f(0))/2)$.
 \end{lem}
  \begin{pf}
  We may assume without loss of generality that
  $f(0)=(\frac{r^{2}-1}{r^{2}+1},0,\ldots,0)$. Let $h$ be a isometry such that $h(0)=f(0).$ The matrix $B\in \mathbf{SU}(n,1)$
  corresponding to $h$ can be written as
  $$B= \left(\begin{array}{ccc}  I_{n-1}& & \\
& \frac{r^{2}+1}{2r}& \frac{r^{2}-1}{2r}\\
&\frac{r^{2}-1}{2r}&
\frac{r^{2}+1}{2r}\end{array}\right).$$
Since $B^{-1}A$ stabilizes the origin, $B^{-1}A\in \mathbf{U}(n)$, then
$$\rm{dist}(A,\mathbf{U}(n))\leq\parallel A-B^{-1}A\parallel\leq  \parallel B^{-1}-Id\parallel \parallel A\parallel.$$
  Also we have
  $$\parallel A\parallel=\parallel BB^{-1}A\parallel \leq  \parallel B^{-1}A\parallel \parallel B\parallel \leq \parallel B\parallel.$$

Simple calculation reveals that
$$\parallel B^{-1}-Id\parallel=r-1$$ and $$\parallel B\parallel=r.$$
\end{pf}
\begin{lem}
Suppose that $f$ is a complex hyperbolic isometry.  Let $A\in \mathbf{SU}(n,1)$ be the matrix corresponding to $f$
and $B\in\mathbf{U}(n)$. Then for each $q\geq 1$
$$\parallel A^{q}-B^{q}\parallel\leq \frac{r^{q}-1}{r-1}\parallel A-B\parallel,$$
where $$r=\exp(\rho(0,f(0))/2).$$
 \end{lem}
\begin{pf}
We recall the identity
$$A^{q}-B^{q}=(A-B)A^{q-1}+B(A-B)A^{q-2}+\ldots B^{q-2}(A-B)A+B^{q-1}(A-B).$$

In the proof of the  Lemma 4.3, we know that $\parallel A \parallel\leq r $. As $B\in\mathbf{U}(n) $, we get $\parallel B^{i}\parallel=1$ for all $i$.
Using the triangle inequality, we have
 $$\parallel A^{q}-B^{q} \parallel \leq \parallel A-B \parallel (r^{q-1}+r^{q-2}+\ldots+r+1)=\frac{r^{q}-1}{r-1}\parallel A-B\parallel.$$

\end{pf}
Now, we prove our main result.
\begin{thm}Let $G\subset \mathbf{PU}(n,1) $ be a discrete torsion-free non-elementary group. Then there exist $o\in \mathbf{H}_{\mathbf{C}}^{n}$ so that
for any $f\in G$
$$\rho\left(o, f(o)\right)\geq \delta, $$ where $\delta=\frac{0.02}{17^{n-1}}$.
\end{thm}

\begin{pf}
First, according to the Lemma 3.3, we know that there exist an $\alpha\in \mathbf{SU}(n,1)$ such that for any $Id\neq C\in \alpha G\alpha^{-1}$
$$\parallel C\parallel \parallel C-Id \parallel\geq \omega\approx 0.3854.$$
We show that the Theorem holds with $o=\alpha^{-1}(0)$. Assume to the contrary that the conclusion does not hold for some $f\in G$.
Denote the corresponding matrix in $\mathbf{SU}(n,1)$ of $f$ be
$\widehat{A}$. Set $A=\alpha \widehat{A}\alpha^{-1}$, $r=\exp(\rho(0,f(0))/2)$. We deduce that $r<e^{\frac{\delta}{2}}$. Furthermore,
there exist $O\in \mathbf{U}(n)$ such that $\parallel A-O \parallel\leq r(r-1)$.

Applying Lemma 4.2 to $O$ with $Q=17$, we then deduce the existence of an elliptic $B\in\mathbf{U}(n) $ of order $q$ such that
$$\parallel B-O \parallel\leq \frac{2\pi}{qQ}, \ \ \parallel O^{q}-B^{q} \parallel=\parallel O^{q}-Id \parallel\leq \frac{2\pi}{Q},$$
 $1\leq q\leq Q^{n-1}.$

By using Lemma 4.4, we have

$$\parallel A^{q}-O^{q} \parallel\leq\frac{r^{q}-1}{r-1}\parallel A-O \parallel\leq r(r^{q}-1).$$

According to the triangle inequality and the fact that $O\in \mathbf{U}(n) $ and the existence of $B$, we can deduce that
$$\parallel A^{q}\parallel\leq(r^{q}-1)r+1$$ and
$$\parallel A^{q}-Id\parallel=\parallel A^{q}-O^{q}+O^{q}-Id\parallel\leq(r^{q}-1)r+\frac{2\pi}{Q}.$$

Thus
$$\parallel A^{q}\parallel\parallel A^{q}-Id\parallel\leq[(r^{q}-1)r+1][(r^{q}-1)r+\frac{2\pi}{Q}].$$

Hence,
$$\parallel A^{q}\parallel\parallel A^{q}-Id\parallel\leq[(e^{\frac{\delta}{2} Q^{n-1}}-1)e^{\frac{\delta}{2}}+1][(e^{\frac{\delta}{2} Q^{n-1}}-1)e^{\frac{\delta}{2}}+\frac{2\pi}{Q}].$$

The assumptions that $Q=17$, $\delta=\frac{0.02}{17^{n-1}}$ and the inequality $n\geq2$ then yield
$$\parallel A^{q}\parallel\parallel A^{q}-Id\parallel\leq[(e^{0.01}-1)e^{0.01/17}+1][(e^{0.01}-1)e^{0.01/17}+\frac{2\pi}{17}]\approx 0.3834\leq 0.3854,$$
which contradicts to Lemma 3.3.
\end{pf}

\begin{cor}
$G\subset \mathbf{PU}(n,1) $ be a discrete non-elementary torsion-free  group of complex hyperbolic isometries of $\mathbf{H}_{\mathbf{C}}^{n}$. Then there is
a Bergman ball of radius
$$r_{n}=\frac{0.01}{17^{n-1}}$$
lying inside every Dirichlet region for $G$.
\end{cor}
Therefore the volume of all complex hyperbolic manifolds $\mathbf{H}_{\mathbf{C}}^{n}/G$ bounded below by the volume of this Bergman ball. We estimate this by
using Lemma 2.1.

\begin{cor}
Let $M$ be a complex hyperbolic $n$-manifolds. Then
$$Vol(M)\geq \frac{4^{n}\sigma_{2n-1}}{2n}\sinh^{2n}(\frac{0.005}{17^{n-1}})$$
where $\sigma_{2n-1}=\frac{2\pi^{n}}{(n-1)!}$ is the Euclidean volume of the unit sphere
$S^{2n-1}\subset \mathbf{C}^{n}$.
\end{cor}

 \vskip 12pt {\bf Acknowledgement}.
 The authors would like to thank the referees for the very helpful comments and
suggestions.
 This research was supported by National Natural Science Foundational of
China (No.10671059). B. Xie also supported by NSF(No.11201134) and 'Young teachers support program' of Hunan University.

\end{document}